\documentclass{amsart}
\usepackage{amssymb, amsmath, tikz}

\def\ZZ{{\mathbb Z}}

\DeclareMathOperator{\Aut}{Aut}
\DeclareMathOperator{\Jac}{Jac}
\DeclareMathOperator{\Pic}{Pic}

\numberwithin{equation}{section}
\newtheorem{theorem}[equation]{Theorem}

\theoremstyle{definition}
\newtheorem{question}[equation]{Question}
\newtheorem{example}[equation]{Example}

\theoremstyle{remark}
\newtheorem{remark}[equation]{Remark}

\begin{document}

\title[Jacobians, Tutte polynomials, and zeta functions]{A note on Jacobians, Tutte polynomials, and two-variable zeta functions of graphs}
\author{Julien Clancy}
\author{Timothy Leake}
\author{Sam Payne}
\address{Department of Mathematics, Yale University, 10 Hillhouse Ave, New Haven CT, 06511 \ \ \tt julien.clancy@yale.edu \ \ timothy.leake@yale.edu \ \ sam.payne@yale.edu}

\begin{abstract}
We address questions posed by Lorenzini about relations between Jacobians, Tutte polynomials, and the Brill--Noether theory of finite graphs, as encoded in his two-variable zeta functions.  In particular, we give examples showing that none of these invariants is determined by the other two.
\end{abstract}

\maketitle

\section{Introduction}

Lorenzini recently introduced the notion of a Riemann--Roch structure on a lattice of rank $n-1$ in $\ZZ^n$ and defined a two-variable zeta function associated to each such structure \cite{Lorenzini12}.  His construction is inspired by earlier work on two-variable zeta functions for number fields and curves over finite fields \cite{Pellikaan96, vanderGeerSchoof00, Deninger03, LagariasRains03}, and includes the Riemann--Roch theory for graphs from \cite{BakerNorine07} as the special case when the lattice is the image of the combinatorial Laplacian of a graph on $n$ vertices.  This note addresses the relationship between his two-variable zeta function in this special case and more classical invariants of graphs, such as the Tutte polynomial and the Jacobian group, which is the torsion part of the cokernel of the combinatorial Laplacian.

Let $G$ be a finite connected graph without loops or multiple edges, and let $\Pic(G)$ be the cokernel of the combinatorial Laplacian of $G$.  Baker and Norine assign a natural degree and rank to each element of $\Pic(G)$, in close analogy with the degree and rank of divisor classes in the Picard group of an algebraic curve \cite{BakerNorine07}, and Lorenzini's two-variable zeta function, which we denote $Z_G(t,u)$, encodes the number of divisor classes in $\Pic(G)$ of each degree and rank.  More precisely,
\[
Z_G(t,u) = \sum_{[D] \in \Pic(G)}  \frac{u^{h(D)} - 1}{u-1} t^{\deg(D)},
\]
where $h(D) = r(D) + 1$, one more than the Baker--Norine rank of $D$; it is the analogue of $h^0$ for a divisor on a smooth algebraic curve.  This zeta function is a rational function and can be expressed as $f_G(t,u) / \big((1-t)(1-tu) \big)$, for some polynomial $f_G$ with integer coefficients.  It also satisfies a functional equation $Z_G(\frac{1}{ut}, u) = (ut^2)^{1-g} Z_G(t,u)$.  Furthermore, $f(1,u)$ is the order of the Jacobian group $\Jac(G)$, the group of divisor classes of degree zero \cite[Proposition~3.10]{Lorenzini12}, which is equal to the number of spanning trees of $G$.

The Tutte polynomial $T_G(x,y)$ also specializes to the order of the Jacobian group, and graphs with the same Tutte polynomial share many other characteristics.  For instance, they have the same number of $k$-colorings for every $k$.  For further details on the relationship between Tutte polynomials, chip-firing, and Jacobians of graphs, see \cite{Gabrielov93, Gabrielov93b, MerinoLopez97, Biggs99}.  Lorenzini asked whether two connected graphs with the same Tutte polynomial must have the same zeta functions or isomorphic Jacobians, and observed that the answers are affirmative for trees \cite[p.~20]{Lorenzini12}.  Our main results are strong negative answers to both of these questions.

\begin{theorem} \label{thm:zeta}
There are pairs of graphs with the same Tutte polynomial and isomorphic Jacobians whose zeta functions are not equal.
\end{theorem}

\begin{theorem} \label{thm:Jac}
There are pairs of graphs with the same Tutte polynomial and the same zeta function whose Jacobians are not isomorphic.
\end{theorem}

\noindent We also give a pair of graphs with the same zeta functions and isomorphic Jacobians whose Tutte polynomials are not equal, so no two of these invariants determines the third.  In each of these pairs, the two graphs also have the same number of vertices and edges.

\begin{remark} \label{rem:method}
Although any two trees on the same number of vertices have the same Tutte polynomial, they also have trivial Jacobians, and it seems to be difficult to construct large classes of pairs of graphs with the same Tutte polynomial and nontrivial Jacobians.  Bollob\'as, Pebody, and Riordan conjectured that the Tutte polynomial (and even its chromatic specialization) is a complete invariant for almost all graphs \cite{BollobasPebodyRiordan00}.

In this project, we pursued two methods for systematically producing graphs with the same Tutte polynomial and nontrivial Jacobians.  One method is exhaustive search, which produced most of the examples in Section~\ref{sec:examples}.  The other method is Tutte's rotor construction, using rotors of order 3, 4, or 5.  A \emph{rotor} is a graph $R$ with an automorphism $\theta$ of the given order $k$ and a vertex $v$ such that $v$, $\theta v$, ..., $\theta^k v$ are distinct.  Tutte's construction takes this rotor together with a map $g$ from the set $\{v, \theta v, \ldots, \theta^k v\}$ to the vertices of another graph $S$ as input.  The output is two new graphs, obtained by gluing $R$ to $S$ in different ways, identifying $\theta^iv$ with either $g(\theta^iv)$ or $g(\theta^{k-i}v)$.  See \cite{Tutte74} for further details and a proof that the resulting two graphs have the same Tutte polynomial.  Note that the same construction with rotors of higher order generally do not produce graphs with the same Tutte polynomial \cite{Foldes78}.  Tutte's rotor construction sometimes produces pairs of graphs whose Jacobians are not isomorphic, as in Example~\ref{ex:rotor}.  Interestingly, however, applying Tutte's construction with his original example of a rotor of order 3 \cite[Figure~2]{Tutte74} has produced pairs of graphs with isomorphic Jacobians in all of our test cases.
\end{remark}

\begin{question} \label{q:rotor}
Does Tutte's construction with his original example of a rotor of order 3 always produce pairs of graphs with isomorphic Jacobians?
\end{question}

\noindent \textbf{Acknowledgments.}  We thank M. Baker, M. Kahle, D. Lorenzini, F. Shokrieh and the referee for helpful conversations and suggestions.  The work of SP is partially supported by NSF DMS--1068689 and NSF CAREER DMS--1149054.

\section{Preliminaries}

Throughout, we consider a finite connected graph $G$, without loops or multiple edges.  Let $v_1, \ldots, v_n$ be the vertices of $G$.  Recall that a divisor on $G$ is a formal sum $D = a_1 v_1 + \cdots + a_n v_n$ with integer coefficients, and the degree of a divisor is the sum of its coefficients
\[
\deg(D) = a_1 + \cdots + a_n.
\]
The combinatorial Laplacian matrix $\Delta(G)$ is the degree matrix minus the adjency matrix of $G$.  Its $i$th diagonal entry is the number of vertices neighboring $v_i$ and its $(i,j)$th off-diagonal entry is $-1$ if $v_i$ is adjacent to $v_j$, and $0$ otherwise.  The combinatorial Laplacian $\Delta(G)$ determines a map from $\ZZ^n$ to the group of divisors whose image is a sublattice of rank $n-1$.  All divisors in this sublattice have degree zero, so the quotient $\Pic(G)$ is graded by degree.  The subgroup consisting of divisor classes of degree zero is called the Jacobian of the graph and denoted $\Jac(G)$.  Since $G$ is connected, $\Jac(G)$ is the torsion subgroup of the cokernel of $\Delta(G)$.  If $G$ has $m$ edges, then the \emph{genus} of $G$ is
\[
g = m - n + 1.
\]

\begin{remark}
This definition of the genus of a graph is the usual one in the literature on tropical geometry and Riemann-Roch, and the terminology reflects a close relation to the genus of certain algebraic curves.  It should not be confused with the minimal genus of a topological surface in which the graph embeds without crossings, which is also called the genus of the graph in the graph theory literature.
\end{remark}

Riemann--Roch theory for graphs, as developed by Baker and Norine in \cite{BakerNorine07}, associates an integer rank $r(D)$ to each divisor $D$ on $G$, analogous to the dimension of a complete linear series on an algebraic curve.  For the purposes of this note, we follow Lorenzini and work with the invariant $h(D) = r(D) + 1$, which is analogous to the dimension of the space of global sections of a line bundle.  It depends only on the class $[D]$ in $\Pic(G)$.  Suppose $D$ has degree $d$.  If $d$ is negative, then $h(D)$ is zero, and if $d > 2g-2$ then $h(D) = d - g + 1$.

Brill--Noether theory, for graphs as for algebraic curves, is concerned with the existence and geometry of divisor classes of given degree and rank.  For a finite graph, one can simply count these classes, and these counts are encoded in Lorenzini's two-variable zeta function
\[
Z_G(t,u) = \sum_{[D] \in \Pic(G)}  \frac{u^{h(D)} - 1}{u-1} t^{\deg(D)}.
\]
Note that two graphs have the same zeta function if and only if they have the same number of divisor classes of each degree and rank.  In particular, any two graphs with the same zeta function have the same number of divisor classes of degree zero, i.e. their Jacobians have the same size.  See \cite{Baker08, tropicalBN, LPP12, Caporaso12b, Len12} for further details on the Brill--Noether theory of graphs.

\section{Examples} \label{sec:examples}

Our first example is a pair of graphs with the same Tutte polynomial and isomorphic Jacobians whose zeta functions are different.

\begin{example}\label{ex:zeta}

Each of the following two graphs is a wedge sum of a triangle with the genus two graph on four vertices; the difference is the vertex of the genus two graph at which the triangle is attached.  Because the Jacobian group of a wedge sum of graphs is the product of the Jacobian groups, and the Tutte polynomial of a wedge sum is the product of the Tutte polynomials, these two graphs have the same Tutte polynomials and isomorphic Jacobians.

\medskip

\begin{center}
\includegraphics{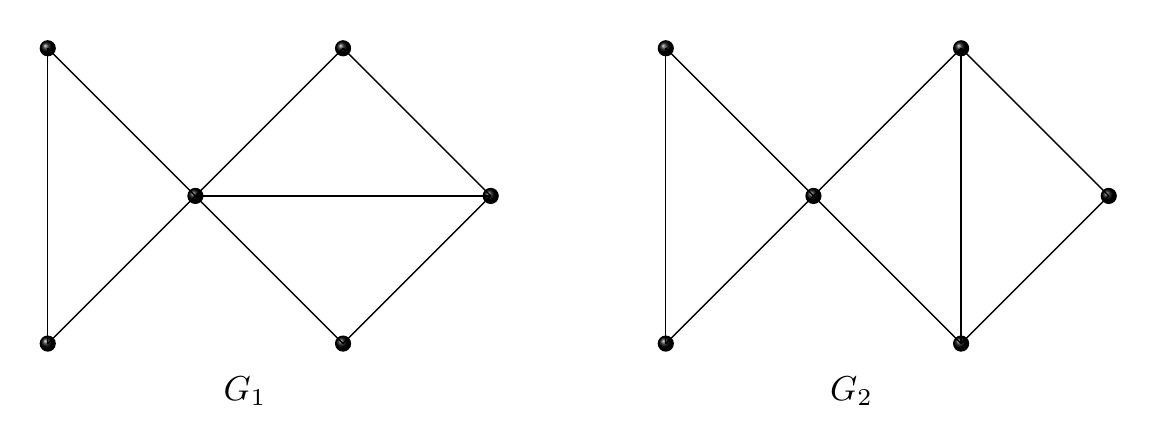}
\end{center}

\noindent One can also compute directly for each graph that the Tutte polynomial is
\[
\left(x+x^2+y\right)\left(x+2x^2+x^3+y+2xy+y^2\right)
\]
and the Jacobian is isomorphic to $\ZZ/24\ZZ$.  However, $G_2$ has a divisor class of degree 2 and rank 1, represented by twice the rightmost vertex, while $G_1$ has no such divisor class.  This can be checked directly, by computing the ranks of each of the 24 divisor classes of degree 2 on $G_1$ and on $G_2$.  Alternatively, one may recall that a graph has a divisor class of degree 2 and rank 1 if and only if it has an involution such that the quotient of the geometric realization by the induced topological involution is a tree \cite{BakerNorine09} and observe that $G_2$ has such an involution, given by a vertical reflection in the figure above, while $G_1$ has no such involution.  It follows that the zeta functions of these two graphs are distinct.  We find that the zeta functions are
\[
Z_{G_1}(t,u) = 1+6t+16t^2+6t^3u+t^4u^2+\frac{24t^3}{(1-t)(1-tu)} \]
\begin{center}
 and
 \end{center}
 \[Z_{G_2}(t,u) = 1+6t+16t^2+t^2u+6t^3u+t^4u^2+\frac{24t^3}{(1-t)(1-tu)}. \qed
\]
\end{example}

\medskip

Our next example is a pair of graphs with the same Tutte polynomial and zeta function whose Jacobians are not isomorphic.

\begin{example} \label{ex:Jac}
Consider the following two graphs of genus $4$ on $8$ vertices.

\begin{center}
\includegraphics{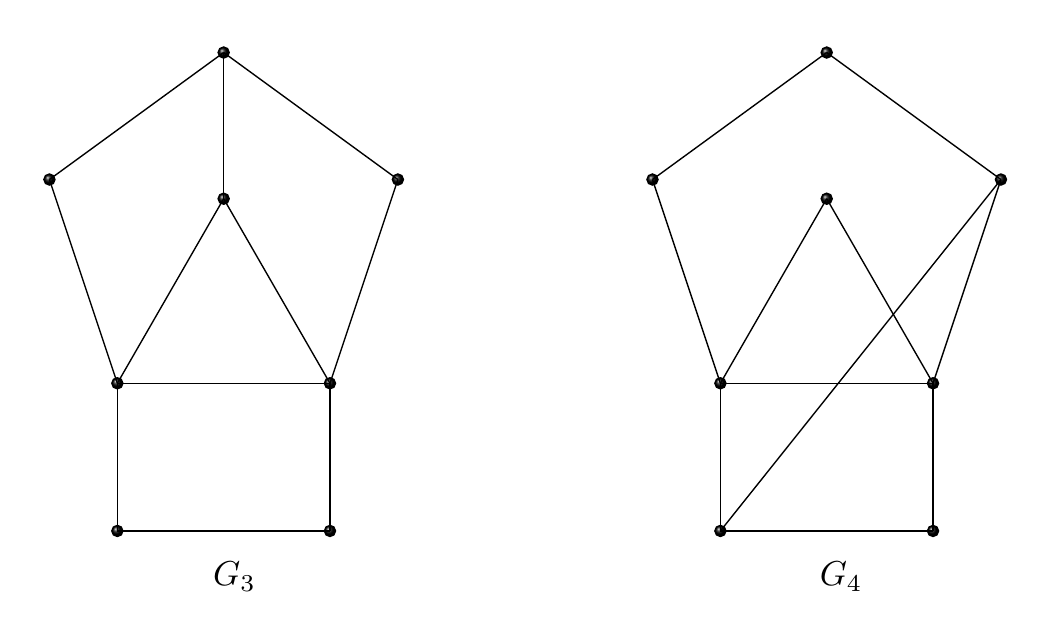}
\end{center}

\noindent Each of these graphs has Tutte polynomial
\begin{eqnarray*}
T(x,y) &=& x^7 + 4x^6 + x^5y + 9x^5 + 6x^4y + 3x^3y^2 + x^2y^3 + 13x^4 + 13x^3y + \\&& + 7x^2y^2 + 3xy^3 + y^4 + 12x^3 + 15x^2y + 9xy^2 + 3y^3 + 7x^2 + \\ && + 9xy + 4y^2 + 2x + 2y,
\end{eqnarray*}
and zeta function
\[
Z(t, u ) = 1+8 t + 31t^2 + 77t^3 + 2t^3u + 31t^4 u+ 8t^5 u^2 + t^6 u^3 + \frac{125t^4}{(1-t)(1-tu)}.
\]
However, their Jacobians are not isomorphic, with
\[
\Jac(G_{3}) \cong \ZZ/5\ZZ \times \ZZ/25\ZZ \mbox{ \ \ and \ \ } \Jac(G_{4}) \cong \ZZ/125\ZZ. \qed
\]
\end{example}

\medskip

\noindent  Examples \ref{ex:zeta} and \ref{ex:Jac} answer Lorenzini's original questions, and prove Theorems~\ref{thm:zeta} and \ref{thm:Jac}, respectively.

\begin{remark}
The graphs in Example~\ref{ex:zeta} have a cut vertex, and those in Example~\ref{ex:Jac} can be disconnected by removing two edges, but there are other (more complicated) examples with higher connectivity whose Tutte polynomials, Jacobians, and zeta functions exhibit similar properties.  For instance, we found a pair of 3-connected graphs of genus 8 on 8 vertices with the same Tutte polynomial and isomorphic Jacobians, whose zeta functions are not equal, as well as a pair of 3-connected graphs of genus 10 on 9 vertices with the same Tutte polynomial whose Jacobians are not isomorphic.

It is not clear what natural graph-theoretic conditions could imply that two graphs with the same Tutte polynomial would also have the same zeta functions, though we did find experimental evidence suggesting that Tutte's rotor construction with certain rotors might produce pairs of graphs with isomorphic Jacobians.  See Remark~\ref{rem:method} and Question~\ref{q:rotor}, above.
\end{remark}

\begin{remark}
Example~\ref{ex:Jac} also gives a negative answer to another question of Lorenzini.  Cori and Rossin proved that planar dual graphs $G$ and $G^*$ have isomorphic Jacobians \cite{CoriRossin00}, and Lorenzini asked \cite[p.~18]{Lorenzini12} whether the existence of this isomorphism follows from the symmetry relating their Tutte polynomials
\[
T_G(x,y) = T_{G^*}(y,x).
\]
The graphs $G_3$ and $G_4$ are planar and have the same Tutte polynomial, so we can choose a planar embedding of $G_4$ to get a dual graph $G_4^*$. The same symmetry holds, $T_{G_3}(x,y) = T_{G_4^*}(y,x)$, but the Jacobians of $G_3$ and $G_4^*$ are not isomorphic.
\end{remark}

We also observe that there are pairs of graphs with the same zeta functions and isomorphic Jacobians whose Tutte polynomials are different, as in the following example.

\begin{example} \label{ex:Tutte}
Consider the following pair of graphs of genus 7 on 7 vertices.

\begin{center}
\includegraphics{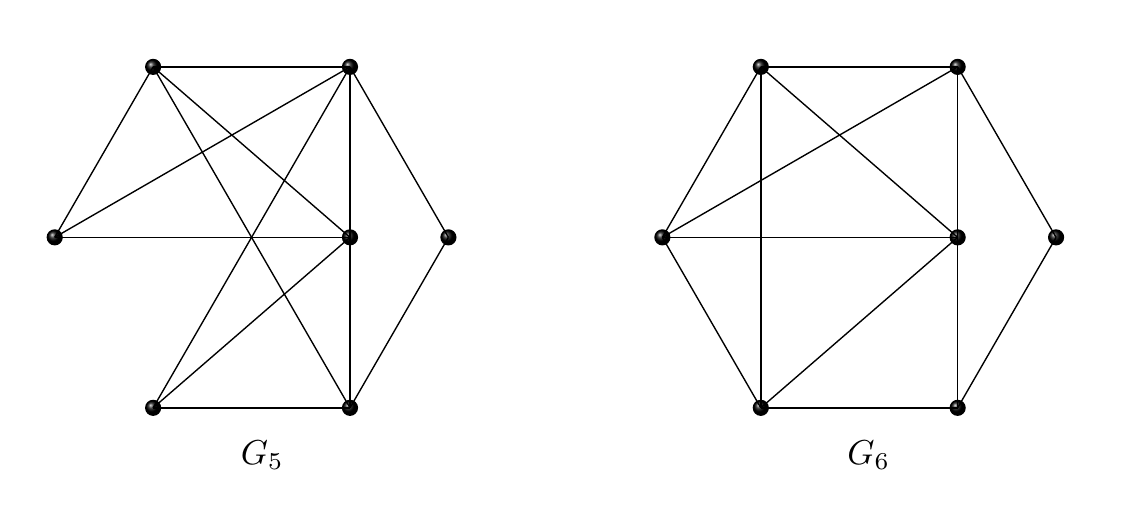}
\end{center}

\noindent We find that $G_5$ and $G_6$ both have Jacobians isomorphic to $\ZZ/545\ZZ$ and zeta functions equal to
\begin{eqnarray*}
Z(t,u) &=& 1+7t+27t^2+75t^3+165t^4+299t^5+449t^6+3t^4u+25t^5u+  \\ && + 105t^6u
+299t^7u+2t^6u^2+25t^7u^2 + 165 t^8 u^2 +3t^8u^3+75t^9u^3+\\ &&+27t^{10}u^4+7t^{11}u^5+t^{12}u^6+\frac{545t^7}{(1-t)(1-tu)}.
\end{eqnarray*}
However, we find that their Tutte polynomials are
\begin{eqnarray*}
T_{G_5}(x,y) & = & 8x+26x^2+33x^3+21x^4+7x^5+x^6+8y  +41xy +60x^2y+ \\ && +34x^3y+7x^4y+23y^2+59xy^2
+43x^2y^2+9x^3y^2+29y^3+\\ && + 44xy^3+16x^2y^3+x^3y^3+23y^4+22xy^4+3x^2y^4+13y^5
+\\ &&+ 7xy^5+5y^6+xy^6+y^7,
\end{eqnarray*}
and
\begin{eqnarray*}
T_{G_6}(x,y) &=& 10x+27x^2+31x^3+20x^4+7x^5+x^6+10y+45xy+55x^2y+\\ &&+32x^3y+8x^4y+28y^2+57xy^2
+38x^2y^2+11x^3y^2+34y^3+\\ &&+ 38xy^3+16x^2y^3 +2x^3y^3+26y^4+17xy^4+5x^2y^4+14y^5
+\\ &&+ 5xy^5+ x^2y^5+5y^6+xy^6+y^7. \qed
\end{eqnarray*}
\end{example}

In the next example, we construct a pair of graphs with the same Tutte polynomial whose Jacobians are not isomorphic using Tutte's rotor construction.

\begin{example} \label{ex:rotor}
We now apply Tutte's rotor construction to the base graph $G_7$ and rotor $G_8$ to construct two graphs of genus 9 on 11 vertices with the same Tutte polynomial.

\vspace{-8pt}

\begin{center}
\includegraphics{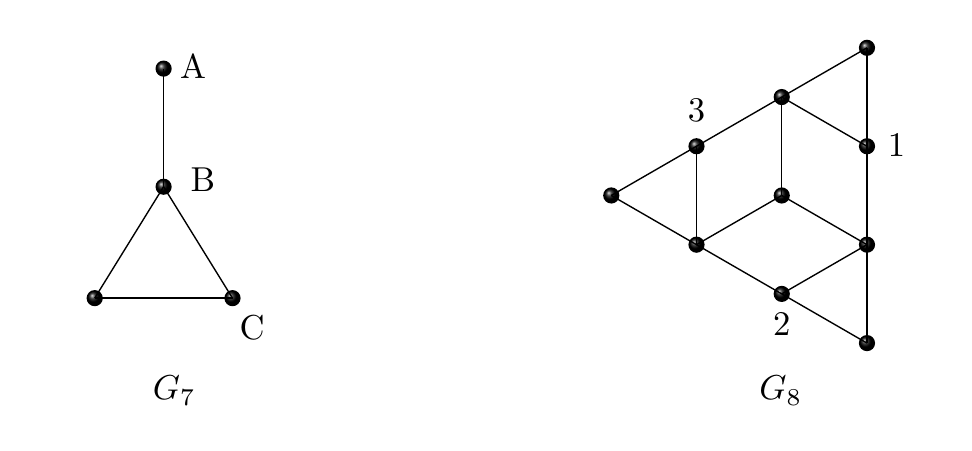}
\end{center}

\vspace{-8pt}

\noindent  The rotor construction involves gluing the vertices $A$, $B$, and $C$ of $G_7$ to the vertices $1$, $2$, and $3$ of $G_8$ in two different ways.  In both cases, we glue $A \mapsto 1$.  For $G_9$, we glue $B \mapsto 2$ and $C \mapsto 3$, whereas for $G_{10}$ we glue $B \mapsto 3$ and $C \mapsto 2$, as shown.

\begin{center}
\includegraphics{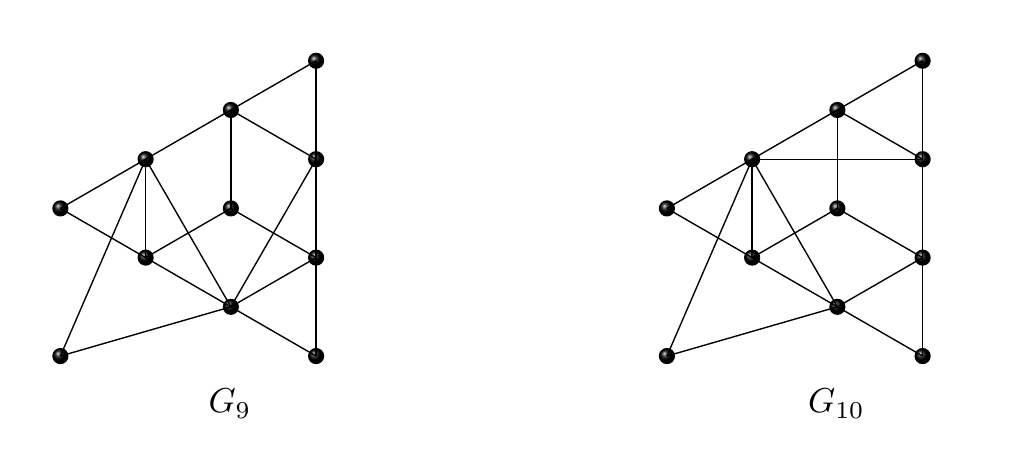}
\end{center}

\noindent By \cite[4.1]{Tutte74}, the graphs $G_9$ and $G_{10}$ have the same Tutte polynomial, which we compute to be
\begin{eqnarray*}
T(x,y) &=& 7 x+47 x^2+139 x^3+239 x^4+266 x^5+202 x^6+107 x^7+39 x^8+9 x^9\\
&&+x^{10}+7 y+72 x y+ 270 x^2 y+525 x^3 y+601 x^4 y+426 x^5 y+188 x^6 y\\
&&+49 x^7 y+6 x^8 y+32 y^2+213 x y^2+553 x^2 y^2+735 x^3 y^2+540 x^4 y^2\\
&&+217 x^5 y^2+43 x^6 y^2+3 x^7 y^2+67 y^3+324 x y^3+597 x^2 y^3+525 x^3 y^3\\
&&+221 x^4 y^3+37 x^5 y^3+x^6 y^3+85 y^4+305 x y^4+392 x^2 y^4+212 x^3 y^4\\
&&+40 x^4 y^4+73 y^5+194 x y^5+167 x^2 y^5+47 x^3 y^5+x^4 y^5+45 y^6+y^9\\
&&+86 x y^6+ 44 x^2 y^6+4 x^3 y^6+20 y^7+25 x y^7+6 x^2 y^7+6 y^8+4 x y^8.\\
\end{eqnarray*}

However, we find that their Jacobians are not isomorphic, with
\[
\Jac(G_9) \cong \ZZ/9065\ZZ \mbox{ \ \  and \ \ } \Jac(G_{10}) \cong \ZZ/1295\ZZ \times \ZZ/7\ZZ.
\]
These two graphs also have distinct zeta functions, with
\begin{eqnarray*}
Z_{G_9} (t,u) &=& 1 + 11t + 62t^2 + 241 t^3 + 723 t^4 + 1757 t^5 + + 3529 t^6+ 5865 t^7 + \\ && + 8009 t^8 + 6 t^5 u  + 86 t^6 u + 589 t^7 u
+ 2385 t^8 u + 5865 t^9 u + \\ &&+ 31 t^8 u^2 + 598 t^9 u^2 + 3529 t^{10} u^2 + 86 t^{10} u^3  + 1757 t^{11} u^3 + \\ && + 6t^{11} u^4
  + 723 t^{12} u^5 + 241 t^{13} u^5  +   62 t^{14} u^6 + 11t^{15} u^7 + \\&& + t^{16} u^8 + \frac{9065t^9}{(1-t)(1-tu)},
 \end{eqnarray*}
\noindent and
\begin{eqnarray*}
Z_{G_{10}} (t, u) &=& 1 + 11t + 62t^2 + 241 t^3 + 723 t^4 + 1757 t^5 + 3529 t^6+ 5865 t^7 + \\&& + 8009 t^8  + 4 t^5 u  + 75 t^6 u + 582 t^7 u
+ 2369 t^8 u + 5865 t^9 u + \\&& + 37 t^8 u^2  + 582 t^9 u^2 + 3529 t^{10} u^2 + 75 t^{10} u^3 + 1757 t^{11} u^3 + \\&& + 4t^{11} u^4
+ 723 t^{12} u^5 + 241 t^{13} u^5 + 62 t^{14} u^6 + 11t^{15} u^7 + \\&& + t^{16} u^8 + \frac{9065t^9}{(1-t)(1-tu)}. \qed
 \end{eqnarray*}

\end{example}

\begin{remark}
By exhaustive search, we find that there are no pairs of graphs on 7 or fewer vertices with the same Tutte polynomial whose Jacobians are not isomorphic.  On 8 vertices, there are 11117 isomorphic classes of connected graphs, and we find only two such pairs.  These are the pairs in Example~\ref{ex:Jac} above, and in Example~\ref{ex:smallest} below.  On 9 vertices, there are 261080 isomorphism classes of connected graphs, but we find only 122 pairs of graphs with the same Tutte polynomial whose Jacobians are not isomorphic.  Some graphs appear in several pairs; these 122 pairs involve 99 different graphs with 33 different Tutte polynomials.
\end{remark}

\begin{example} \label{ex:smallest}
As mentioned above, there are only two pairs of graphs on 8 vertices with the same Tutte polynomial whose Jacobians are not isomorphic.  One such pair is given in Example~\ref{ex:Jac}.  The other is the following pair of graphs of genus 5, which we include for completeness.

\begin{center}
\includegraphics{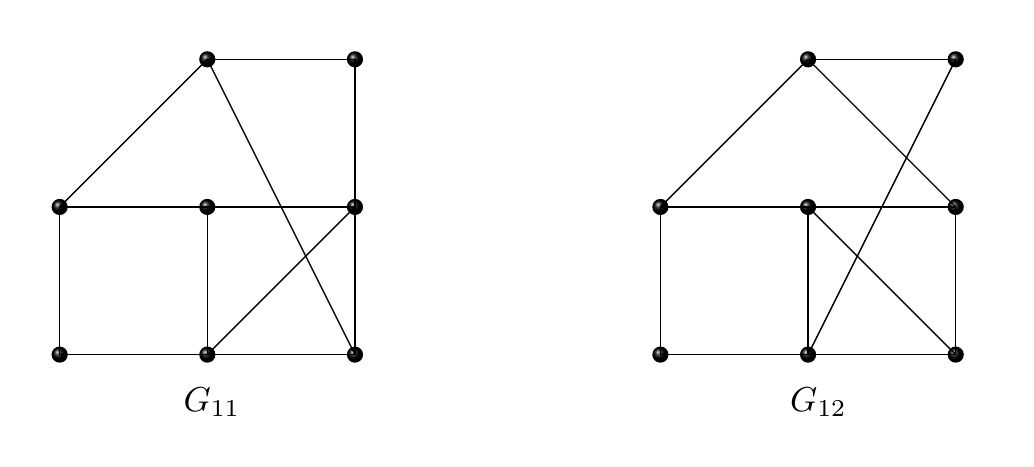}
\end{center}

\noindent Each of these graphs has Tutte polynomial
\begin{eqnarray*}
T(x,y) &=& 5x+18x^2+27x^3+23x^4+13x^5+5x^6+x^7+5y+25xy +38x^2y +   \\ && +25x^3y+9x^4y+2x^5y + 12y^2 +28xy^2+18x^2y^2+5x^3y^2+x^4y^2+ \\ &&  +11y^3 + 12xy^3+3x^2y^3+5y^4+2xy^4+y^5,
\end{eqnarray*}

\noindent and zeta function

\begin{eqnarray*}
Z(t,u) &=&  1+8t+34t^2+98t^3+202t^4+13t^4u+98t^5u+34t^6u^2 + \\ && + 8t^7u^3  +t^8u^4+\frac{294t^5}{(1-t)(1-tu)}.
\end{eqnarray*}
However, we find that their Jacobians are not isomorphic, with
\[
\Jac(G_{11}) \cong \ZZ/42\ZZ\times \ZZ/7\ZZ \mbox{ \ \ and \ \ } \Jac(G_{12}) \cong \ZZ/294\ZZ.  \qed
\]

\end{example}

\begin{remark}
Gim\'enez and Merino independently found a pair of graphs with the same Tutte polynomial whose Jacobians are not isomorphic \cite{GimenezMarino02}.  Their example consists of two planar graphs of genus 10 on 12 vertices whose zeta functions are not equal.
\end{remark}

\section{Jacobians of random graphs}

Jacobians of graphs are frequently cyclic, as has been observed by Lorenzini and others.  Perturbing a graph with non-cyclic Jacobian slightly, by subdividing an edge, tends to produce graphs with cyclic Jacobians.  See, for instance, the computation of cyclic Jacobians for modified wheel graphs in \cite{Biggs07}.  This phenomenon was observed again by Robeva in 2008, in computations related to the tropical proof of the Brill--Noether Theorem \cite{tropicalBN}, and we encountered it once more through extensive computations exploring the questions addressed in this paper.  

In the course of this investigation, we observed a somewhat more precise structure.  The Jacobian group of a graph comes with a canonical duality pairing \cite{Shokrieh10}, and a finite abelian group with duality pairing $(G, \langle \ , \ \rangle)$ seems to appear with frequency proportional to
\[
\frac{1}{|G| \  | \Aut(G, \langle \ , \ \rangle )|}.
\]
Here, $\Aut(G, \langle \ , \ \rangle )$ denotes the subgroup of the automorphism group of $G$ that preserves the pairing.  This experimentally observed variation on the Cohen--Lenstra heuristic \cite{CohenLenstra84} should explain the prevalence of cyclic Jacobians.  We find that the Jacobian of a random graph is cyclic with probability slightly greater than $.79$.  Further details on this heuristic and our experiments with Jacobians of random graphs are presented in \cite{RandomJacobians}, some of our conjectures based on this heuristic are proved in \cite{Wood14}, and connections to symmetric function theory and Hall-Littlewood polynomials are studied in \cite{Fulman14}.

\newcommand{\etalchar}[1]{$^{#1}$}
\providecommand{\bysame}{\leavevmode\hbox to3em{\hrulefill}\thinspace}
\providecommand{\MR}{\relax\ifhmode\unskip\space\fi MR }
\providecommand{\MRhref}[2]{%
  \href{http://www.ams.org/mathscinet-getitem?mr=#1}{#2}
}
\providecommand{\href}[2]{#2}

\end{document}